\documentclass[12pt]{amsart}
\usepackage{amssymb}
\usepackage[pagewise]{lineno}
\newtheorem{theorem}{Theorem}[section]
\newtheorem{lemma}[theorem]{Lemma}

\newtheorem{corollary}[theorem]{Corollary}
\newtheorem{proposition}[theorem]{Proposition}

\newtheorem{lem-def}[theorem]{Lemma-Definition}

\DeclareRobustCommand\longtwoheadrightarrow
{\relbar\joinrel\twoheadrightarrow}

\newcommand{\hooklongrightarrow}{\lhook\joinrel\longrightarrow}
\renewenvironment{proof}{{\bfseries Proof.}}{\qed}
\newcommand{\N}{\mathbb N}
\newcommand{\Z}{\mathbb Z}
\newcommand{\Q}{\mathbb Q}

\def\op{\operatorname}

\def\aa{\mathcal{A}}

\def\al{\alpha}

\def\ars#1{\renewcommand\arraystretch{#1}}

\def\aut{\op{Aut}}

\def\bad{B(a,\delta)}
\def\bb{{\mathcal B}}

\def\be{\beta}

\def\cc{\mathcal{C}}

\def\cfa{\left(\ri\right)_{i\in A}}
\def\cfb{\left(\zj\right)_{j\in B}}

\def\com{{\op{com}}}

\def\cl#1{\left[\,#1\,\right]_\mu}

\def\dd{D_\delta}
\def\ddp{D_{\delta'}}
\def\defn{\nn{\bf Definition. }}

\def\dg{\op{deg}}

\def\diso{\lower.4ex\hbox{$\downarrow$}\raise.4ex\hbox{\mbox{\scriptsize
$\wr$}}}

\def\dta{\delta}

\def\e{\medskip}

\def\ep#1{\exp(\Pi i#1)}
\def\ep{\epsilon}
\def\epm{\epsilon_\mu}
\def\epn{\epsilon_\nu}

\def\fin{{\op{fin}}}
\def\fmp{F_{\mu,\phi}}

\def\g{\Gamma}
\def\ga{\gamma}

\def\gal{\op{Gal}}

\def\gen#1{\big\langle\, {#1} \,\big\rangle}

\def\gg{\mathcal{G}}

\def\ggm{\mathcal{G}_\mu}

\def\ggo{\mathcal{G}_\om}

\def\gm{\g_\mu}

\def\gn{\g_\nu}

\def\gq{\g_\Q}

\def\gv{\Gamma_v}

\def\hk{\hookrightarrow}

\def\imp{\ \Longrightarrow\ }

\def\inm{\op{in}_\mu}
\def\inn{\op{in}}

\def\ino{\op{in}_\om}

\def\irr{\op{Irr}}
\def\ism{\lower.3ex\hbox{\ars{.08}$\begin{array}{c}\,\to\\\mbox{\tiny $\sim\,$}\end{array}$}}
\def\iso{\ \lower.3ex\hbox{\ars{.08}$\begin{array}{c}\lra\\\mbox{\tiny $\sim\,$}\end{array}$}\ }

\def\kb{\overline{K}}
\def\kbx{\overline{K}[x]}
\def\kh{K^h}
\def\khx{K^h[x]}

\def\kp{\op{KP}}

\def\kpi{\op{KP}_\infty}
\def\kpm{\op{KP}(\mu)}

\def\kpn{\op{KP}(\nu)}
\def\kpo{\op{KP}(\om)}

\def\ks{K^{\op{sep}}}

\def\kx{K[x]}

\def\La{\Lambda}

\def\lfin{\ll_\fin}
\def\lg{l\raise.6ex\hbox to.2em{\hss.\hss}l}
\def\li{\ll_\infty}
\def\ll{\mathcal{L}}
\def\lra{\,\longrightarrow\,}

\def\ml{\op{mult}}

\def\mmu{\mid_\mu}

\def\muh{\mu^h}

\def\nb{\bar{\nu}}

\def\nn{\noindent}

\def\nuh{\nu^h}

\def\om{\omega}

\def\orb{\hbox to  .3em{$\backslash$}\backslash}

\def\p{\mathfrak{p}}

\def\pb#1{\partial_b#1}

\def\ppa{\mathcal{P}_{\al}}

\def\pset{\mathcal{P}}

\def\res{\operatorname{res}}

\def\rha{\rho_\aa}

\def\rhc{\rho_\cc}
\def\ri{\rho_i}

\def\rj{\rho_j}

\def\scut{\sim_{\op{qcut}}}
\def\sg{\sigma}

\def\sii{\ \Longleftrightarrow\ }

\def\smu{\sim_\mu}
\def\snu{\sim_\nu}

\def\sp{\op{Spec}}

\def\supp{\op{supp}}

\def\tinn{\ttt^{\op{inn}}}
\def\tinnk{\ttt^{\op{inn}}(K)}

\def\tk{\ttt(K)}
\def\tkb{\ttt(\kb)}
\def\tkh{\ttt(\kh)}

\def\ttt{\mathcal{T}}
\def\tv{\mathcal{T}_v}
\def\tvh{\mathcal{T}_{\vh}}
\def\tvinn{\tv^{\op{inn}}}
\def\tws{\mathcal{T}^{\op{ws}}}
\def\twsk{\mathcal{T}^{\op{ws}}(K)}
\def\twskh{\mathcal{T}^{\op{ws}}(\kh)}
\def\ty{\mathbf{t}}

\def\vb{\bar{v}}
\def\vh{v^h}

\def\zj{\zeta_j}

\newcounter{cs}
\stepcounter{cs}
\newcommand{\casos}{\begin{itemize}}
\newcommand{\fcasos}{\end{itemize}\setcounter{cs}{1}}

\newfont{\tit}{cmr12 scaled \magstep3}

\title[Rigidity of valuative trees]{Rigidity of valuative trees under henselization}

\makeatletter
\@namedef{subjclassname@2010}{%
  \textup{2010} Mathematics Subject Classification}
%, 12J10}

\author[Nart]{Enric Nart}
\address{Departament de Matem\`{a}tiques,         Universitat Aut\`{o}noma de Barcelona,         Edifici C, E-08193 Bellaterra, Barcelona, Catalonia}
\email{nart@mat.uab.cat}
\thanks{Partially supported by grant  PID2020-116542GB-I00  from the Spanish Research Agency}
\date{\today}

\begin{document}
%MTM2016-75980-P, MTM2015-69135-P

%\pagestyle{empty}

%No assumption is made on the rank of the valuations.
\begin{abstract}
Let  $(K,v)$ be a valued field and let $(\kh,\vh)$ be the henselization determined by the choice of an extension of $v$ to an algebraic closure of $K$.  Consider an embedding $v(K^*)\hk\La$ of the value group into a divisible ordered abelian group. Let $\ttt(K,\La)$, $\ttt(\kh,\La)$ be the trees formed by all $\La$-valued extensions of $v$, $\vh$ to the polynomial rings $\kx$, $\khx$, respectively. We show that the natural restriction mapping $\ttt(\kh,\La)\to\ttt(K,\La)$ is an isomorphism of posets. 

As a consequence, the restriction mapping $\ttt_{\vh}\to\ttt_v$ is an isomorphism of posets too, where $\ttt_v$, $\ttt_{\vh}$ are the trees whose nodes are the equivalence classes of valuations on $\kx$, $\khx$ whose restriction to $K$, $\kh$ are equivalent to $v$, $\vh$, respectively.  

%For any given extension $\mu$ of $v$ to , there exists a unique common extension of $\vh,\mu$ to $\khx$. More precisely, for we find an isomorphism of posets 

\end{abstract}

\maketitle

%\begin{center}\sl Preliminary version \end{center}

%\tableofcontents

\section*{Introduction}
Let $(K,v)$ be a valued field, $\kb$ an algebraic closure of $K$, and $\vb$ an extension of $v$ to $\kb$. Let $(\kh,\vh)$ be the henselization of $(K,v)$ determined by the choice of $\vb$.

The pioneering work of MacLane \cite{mcla,mclb} describing the extensions of $v$ to the polynomial ring $\kx$, in the discrete rank-one case, was  generalized to arbitrary valued fields, independently by Vaqui\'e \cite{Vaq0,Vaq} and Herrera, Mahboub, Olalla and Spivakovsky \cite{hos,hmos}.

In-between of these achievements, Alexandru, Popescu and Zaharescu studied the extensions of $v$ to $\kx$ as restrictions of valuations on $\kbx$ extending $\vb$ \cite{APZ0,APZ}, whose description is much easier because most of them are monomial valuations. A similar approach was followed by Kuhlmann \cite{Kuhl}.
Besides these contributions, we find in the literature an extensive analysis of the restriction and lifting of valuations from $\kbx$ to $\kx$ and viceversa \cite{MMS,N2019,PP,Vaq3}.

In this paper, we use these techniques to prove that for
every extension $\mu$ of $v$ to $\kx$ there is a unique common extension of $\mu$ and $\vh$ to $\khx$.

More precisely, in Section \ref{secKh} we prove the following result.\e

\noindent{\bf Theorem A. }{\it Let $\La$ be a divisible ordered abelian group extending $v(K^*)$. Consider the trees $\ttt(K,\La)$, $\ttt(\kh,\La)$ of all $\La$-valued extensions of $v$, $\vh$ to $\kx$, $\khx$, respectively. Then, the natural restriction mapping $\ttt(\kh,\La)\to\ttt(K,\La)$ is an isomorphism of posets.} \e

As a consequence, we derive in Section \ref{secEquiv} a similar result for the trees whose nodes are equivalence classes of  
valuations, without fixing any ordered group containing all value groups.

Two valuations $\mu$, $\nu$ on $\kx$ are \emph{equivalent} if there is an order preserving isomorphism between their value groups $\iota\colon \gm \ism\gn$, such that $\nu=\iota\circ \mu$. \e

\noindent{\bf Theorem B. }{\it  Consider the trees $\ttt_v$, $\ttt_{\vh}$ whose nodes are the equivalence classes of valuations on $\kx$, $\khx$ whose restriction to $K$, $\kh$ are equivalent to $v$, $\vh$, respectively.  
 Then, the restriction mapping $\ttt_{\vh}\to\ttt_v$ is an isomorphism of posets.} \e

Sections \ref{secVals} and \ref{secKb} contain some basic facts on key polynomials, valuative trees and lifting of valuations fom $\kx$ to $\kbx$.

%\newpage

\section{Valuative trees}\label{secVals}

In this section, we recall some well-known results on valuations on a polynomial ring and valuative trees.

Let $(K,v)$ be a valued field and  $\,\g=\gv=v(K^*)\,$ its value group.
We fix throughout the paper an embedding $\g\hk\Lambda$ of ordered abelian groups. We write simply $\Lambda\infty$ instead of $\Lambda\cup\{\infty\}$.

Let $\ttt=\ttt((K,v),\La)$ be the set of all valuations $$\mu\colon \kx\lra\La\infty,$$ whose restriction to $K$ is $v$. 
For any $\mu\in\ttt$, the \emph{support} of $\mu$ is the prime ideal $$\p=\supp(\mu)=\mu^{-1}(\infty)\in\sp(\kx).$$ 

The valuation $\mu$ induces in a natural way a valuation $\bar{\mu}$ on the residue field $K(\p)$, the field of fractions of $\kx/\p$. 
Note that $K(0)=K(x)$, while $K(\p)$ is a simple finite extension of $K$ if $\p\ne0$.

For any given field $L$ we shall denote by $\irr(L)$ the set of monic irreducible polynomials in $L[x]$.

\subsection{Key polynomials and augmentations}\label{subsecKP}
The value group of $\mu$ is the subgroup $\gm\subset\La$ generated by $\mu\left(\kx\setminus\p\right)$. 
We say that $\mu$ is \emph{commensurable} (over $v$) if $\g_\mu/\g$ is a torsion group. In this case, there is a canonical embedding of $\gm$ into the divisible hull of $\g$:
$$\gm\hooklongrightarrow \gq:=\g\otimes_\Z\Q.$$ 

For any $\alpha\in\g_\mu$, consider the abelian groups:
$$
\ppa=\{g\in \kx\mid \mu(g)\ge \alpha\}\supset
\ppa^+=\{g\in \kx\mid \mu(g)> \alpha\}.
$$    
The \emph{graded algebra of $\mu$} is the integral domain:
$$
\ggm=\bigoplus\nolimits_{\alpha\in\g_\mu}\ppa/\ppa^+.
$$

We say that $\ggm$ is \emph{simple} if all nonzero homogeneous elements are units.

There is an \emph{initial term} mapping $\inm\colon \kx\to \ggm$, given by $\inm\p=0$ and 
$$
\inm g= g+\pset_{\mu(g)}^+\in\pset_{\mu(g)}/\pset_{\mu(g)}^+, \qquad\mbox{if }\ g\in \kx\setminus\p.
$$

The following definitions translate properties of the action of  $\mu$ on $\kx$ into algebraic relationships in the graded algebra $\ggm$.\e

\defn Let $g,\,h\in \kx$.

We say that $g,h$ are \emph{$\mu$-equivalent}, and we write $g\smu h$, if $\inm g=\inm h$. 

We say that $g$ is \emph{$\mu$-divisible} by $h$, and we write $h\mmu g$, if $\inm h\mid \inm g$ in $\ggm$.

We say that $g$ is $\mu$-irreducible if $(\inm g)\ggm$ is a nonzero prime ideal. 

We say that $g$ is $\mu$-minimal if $g\nmid_\mu f$ for all nonzero $f\in \kx$ with $\deg(f)<\deg(g)$.\e

%The $\mu$-minimality condition admits a relevant characterization.

\defn A  \emph{(MacLane-Vaqui\'e) key polynomial} for $\mu$ is a monic polynomial in $\kx$ which is simultaneously  $\mu$-minimal and $\mu$-irreducible. It is necessarily irreducible in $\kx$. 
The set of key polynomials for $\mu$ is denoted $\kpm\subset\irr(K)$. \e

The set $\ttt$ of all $\La$-valued valuations on $\kx$ extending $v$, admits a partial ordering
$$\mu\le\nu\ \sii\ \mu(f)\le \nu(f), \quad\forall\,f\in\kx.$$
As usual, we write $\mu<\nu$ to indicate that $\mu\le\nu$ and $\mu\ne\nu$.

%If $\mu\le\nu$, there is a canonical homomorphism of graded $\gg_v$-algebras:$$\ggm\lra\ggn,\qquad \inm f\longmapsto\begin{cases}\inu f,& \mbox{ if }\mu(f)=\nu(f),\\ 0,& \mbox{  if }\mu(f)<\nu(f).\end{cases}$$

This poset $\ttt$ has the structure of a tree. By this, we simply mean that the intervals 
$$
(-\infty,\mu\,]:=\left\{\rho\in\ttt\mid \rho\le\mu\right\}
$$
are totally ordered \cite[Thm. 2.4]{MLV}. 

A node $\mu\in\ttt$ is a \emph{leaf} if it  is a maximal element with respect to the ordering $\le$. Otherwise, we say that $\mu$ is an \emph{inner node}.

The following result combines \cite[Thm. 2.3]{MLV} and \cite[Thm. 4.4]{KP}.

\begin{theorem}\label{maximal}
	A node $\mu\in\ttt$ is a leaf if and only if $\,\kpm=\emptyset$. In this case, the graded algebra $\ggm$ is simple.
\end{theorem}

All valuations with nontrivial support are leaves of $\ttt$. We call them \emph{finite leaves}.
The leaves of $\ttt$ having trivial support are called \emph{infinite leaves}. We denote 
$$
\ttt=\tinn\,\sqcup\,\lfin\,\sqcup\,\li
$$
the subsets of inner nodes, finite leaves, and infinite leaves, respectively.\e

\defn
Take $\mu \in \tinn\sqcup \lfin$. If $\mu$ is an inner node, let $\phi$ be any key polynomial for $\mu$ of minimal degree. If $\mu$ is a finite leaf, let $\phi$ be the monic irreducible generator of $\supp(\mu)$. 
Then, we define the \emph{degree} of $\mu$ as $\deg(\mu)=\deg(\phi)$. \e

A \emph{tangent direction} of $\mu\in\tinn$ is a $\mu$-equivalence class $\cl{\phi}\subset\kpm$ determined by all key polynomials having the same initial term in $\ggm$.

\begin{lemma}\cite[Prop. 2.2]{MLV}\label{propertiesTMN}
	Let $\nu<\mu$ be two nodes in $\ttt$. The set $\ty(\nu,\mu)$  of all monic polynomials $\phi\in\kx$ of minimal degree satisfying $\nu(\phi)<\mu(\phi)$ is a tangent direction of $\nu$. %: \ $\tmn\in\tdm$.
Moreover, for any $\phi\in \ty(\nu,\mu)$ and all $f\in\kx$, the equality $\nu(f)=\mu(f)$ holds if and only if $\phi\nmid_{\nu}f$.   
\end{lemma}

For any $\phi\in\kpn$ and any $\ga\in\Lambda\infty$ such that $\nu(\phi)<\ga$, we may construct the \emph{ordinary augmented valuation}  $\mu=[\nu;\,\phi,\ga]\in\ttt$,
defined in terms of $\phi$-expansions as
$$
f=\sum\nolimits_{0\le s}a_s\phi^s,\ \deg(a_s)<\deg(\phi)\ \imp\ \mu(f)=\min\{\nu(a_s)+s\ga\mid 0\le s\},
$$
Note that $\mu(\phi)=\ga$, $\nu<\mu$ and  $\ty(\nu,\mu)=[\phi]_\nu$. %For any $f\in\kx$ we have$$\nu(f)=\mu(f)\sii \phi\nmid_\nu f.$$
If $\ga<\infty$, then $\phi$ is a key polynomial for $\mu$ of minimal degree \cite[Cor. 7.3]{KP}. 

Consider a totally ordered family 
$
\cc=\cfa\subset\ttt,
$
parametrized by a totally ordered set $A$ of indices such that  the mapping
$A\to\cc$ determined by $i\to \ri$
is an isomorphism of totally ordered sets.

If $\deg(\rho_{i})$ is stable for all sufficiently large $i\in A$, we say that $\cc$ has \emph{stable degree}, and we denote this stable degree by $\dg(\cc)$. 

We say that $f\in\kx$ is \emph{$\cc$-stable} if for some index $i\in A$, it satisfies
$$\ri(f)=\rj(f), \quad \mbox{ for all }\ j>i.$$
We obtain a \emph{stability function} $\rhc$, defined
on the set of all $\cc$-stable polynomials by $\rhc(f)=\max\{\ri(f)\mid i\in A\}$.

A \emph{limit key polynomial} for $\cc$ is a monic $\cc$-unstable polynomial of minimal degree. Let $\kpi(\cc)$ be the set of all these limit key polynomials. Since the product of stable polynomials is stable, all limit key polynomials are irreducible in $ \kx$.\e
%The \emph{limit degree} of $\cc$, denoted  $\dgi(\cc)$,  is the degree of any limit key polynomial. \e

\defn
We say that $\cc$ is an \emph{esssential continuous family} of valuations in $\ttt$ if it satisfies the following conditions:
\begin{enumerate}
\item[(a)] It contains no maximal element and has stable degree.
\item[(b)] It admits limit key polynomials whose degree is greater than $\dg(\cc)$.
\end{enumerate}\e

Take any limit key polynomial $\phi\in\kpi\left(\cc\right)$, and any $\ga\in\La\infty$ such that $\ri(\phi)<\ga$ for all $i\in A$.
We define the \emph{limit augmentation}   $\mu=[\cc;\,\phi,\ga]\in\ttt$ as the following mapping, defined in terms of $\phi$-expansions: 
$$ f=\sum\nolimits_{0\le s}a_s\phi^s\imp \mu(f)=\min\{\rhc(a_s)+s\ga\mid 0\le s\}.
$$
Since $\deg(a_s)<\deg(\phi)$, all coefficients $a_s$ are $\cc$-stable. 
 
Note that $\mu(\phi)=\ga$ and $\ri<\mu$ for all $i\in A$.
If $\ga<\infty$, then $\phi$ is a key polynomial for $\mu$ of minimal degree \cite[Cor. 7.13]{KP}.

\subsection{Well-specified valuations}\label{subsecWell}

A valuation $\mu\in\ttt$ is said to be \emph{well-specified} if it is not an infinite leaf. That is, $\mu$ is either an inner node or a finite leaf.

We denote the subtree of $\ttt$ formed by all well-specified valuations by 
$$
\tws=\tinn\sqcup\lfin.
$$

By the main result of MacLane-Vaqui\'e, the well-specified valuations are those which can be  attained after a finite number of augmentations starting with a monomial (degree one) valuation  \cite{Vaq,MLV}. In particular, they can be obtained by a single augmentation step, either ordinary or limit. In other words, there exists either $\nu\in\tinn$ and $\phi\in\kpn$, or an essential continuous family $\cc$ and $\phi\in\kpi(\cc)$, such that   
$$
\mu=[\nu;\,\phi,\ga],\quad \mbox{or}\quad \mu=[\cc;\,\phi,\ga],
$$  
for some $\ga\in\La$.\e

\noindent{\bf Definition. }In this situation, we say that the polynomial $\phi\in\irr(K)$ \emph{defines} the valuation $\mu$.\e

Caution: this is not an intrinsic property of $\phi$. Any $\phi\in\irr(K)$ may define many different well-specified valuations. 

\begin{lemma}\label{defining}
Let $\mu\in\ttt$ be a well-specified valuation. A monic irreducible polynomial $\phi\in\irr(K)$ defines $\mu$ if and only if one of the following conditions holds:
\begin{enumerate}
\item[(i)] $\mu\in\lfin$ \;and $\,\supp(\mu)=\phi\kx$. 
\item[(ii)] $\mu\in\tinn$ and $\,\phi$ is a key polynomial of minimal degree for $\mu$.
\end{enumerate}  	
\end{lemma}

\begin{proof}
By the theorem of Maclane--Vaqui\'e, we have $\mu=[\nu;\,\phi,\ga]$ or $\mu=[\cc;\,\phi,\ga]$.

If $\mu$ is a finite leaf, we must have $\mu(\phi)=\ga=\infty$, because otherwise $\mu$ would have trivial support. Since $\phi$ is monic and irreducible, we have $\supp(\mu)=\phi\kx$. On the other hand, this property characterizes $\phi$ as the unique monic irreducible polynomial which generates the prime ideal  $\supp(\mu)$ in $\kx$.

Suppose that $\mu$ is an inner node, or equivalently, $\kpm\ne\emptyset$. By \cite[Cor. 7.3, 7.13]{KP} $\phi$ is a key polynomial for $\mu$ of minimal degree.

Conversely, let $\varphi\in\kpm$ be another key polynomial of minimal degree. We may write $\varphi=\phi+a$, for some $a\in\kx$ with $\deg(a)<\deg(\phi)$. Since $\inn_\mu \varphi$ is an homogeneous prime element in $\ggm$, we must have
$\mu(a)\ge \ga=\mu(\phi)$,
because otherwise $\varphi\smu a$ would be a unit in $\ggm$ \cite[Prop. 3.5]{KP}. 

If $\mu=[\nu;\,\phi,\ga]$, we deduce $\nu(a)=\mu(a)\ge\ga>\nu(\phi)$. Hence, $\phi\snu \varphi$. By \cite[Lem. 2.5]{KP}, $\varphi$ is a key polynomial for $\nu$. Then, \cite[Lem. 2.8]{MLV} shows that $\mu=[\nu;\,\varphi,\ga]$, so that $\varphi$ defines $\mu$.

Now, suppose that $\mu=[\cc;\,\phi,\ga]$. Since $\deg(a)<\deg(\phi)$, the polynomial $a$ is $\cc$-stable. Its stable value $\rho_\cc(a)$ satisfies
$$\rho_\cc(a)=\mu(a)\ge\ga=\mu(\phi)>\ri(\phi),\quad \forall i\in A.$$ By \cite[Lem. 4.8]{VT}, $\varphi$ is a limit key polynomial for $\cc$. Then, \cite[Lem. 3.7]{MLV} shows that $\mu=[\cc;\,\varphi,\ga]$, so that $\varphi$ defines $\mu$.
\end{proof}\e

\section{Lifting valuations from $\kx$ to $\kbx$}\label{secKb}
We fix an algebraic closure $\kb$ of $K$, and an extension $\vb$ of our base valuation $v$ to $\kb$. This determines a henselization $(\kh,\vh)$ of $(K,v)$. If $\ks$ is the separable closure of $K$ in $\kb$, the field $\kh$ is the fixed field of the decomposition group 
$$
D_{\vb}=\left\{\sg\in \gal(\ks/K)\mid \vb\circ\sg=\vb\right\}.
$$
The valuation $\vh$ is just the restriction of $\vb$ to $\kh$.
We have a chain of fields
$$
K\subset\kh\subset\ks\subset\kb,
$$
and since $\kb/\ks$ is purely inseparable, we have
$$
\vb\circ\sg=\vb,\quad\mbox{ for all }\sg\in\aut(\kb/\kh). 
$$ 
Thus, the valuation $\vh$ has a unique extension to $\kb$.

From now on, we assume that $\La$ is a divisible group. There are a unique isomorphism $\g_{\vb}\ism\gq$ and a unique embedding $\gq\hk\La$ compatible with the embeddings $\g\hk\g_{\vb}$ and $\g\hk\La$. Hence, up to identifying all groups with their images in $\La$, we have a chain of inclusions  $\gv=\g\subset\g_{\vb}=\gq\subset\La$.

Consider the trees of $\La$-valued valuations on $\kx$, $\khx$, $\kbx$, whose restriction to $K$, $\kh$, $\kb$ is $v$, $\vh$, $\vb$, respectively: 
$$
\tk=\ttt\left((K,v),\La\right),\quad
\tkh=\ttt\left((\kh,\vh),\La\right),\quad
\tkb=\ttt\left((\kb,\vb),\La\right).
$$
Each tree has its corresponding subsets of inner nodes, finite leaves and infinite leaves
$$
\ttt(L)=\tinn(L)\,\sqcup\,\lfin(L)\,\sqcup\,\li(L),
$$
for $L=K$, $\kh$, $\kb$. As shown in \cite[Lem. 3.2]{Kuhl}, the natural restriction mappings between these trees are onto:
$$
\tkb\stackrel{\res}\longtwoheadrightarrow\tkh\stackrel{\res}\longtwoheadrightarrow\tk.
$$
Also, they respect the subsets of inner nodes, finite leaves and infinite leaves.
Indeed, take any $\mu$ in $\tkb$ or $\tkh$. It is obvious that $\mu$ is a finite leaf if and only if $\res(\mu)$ is a finite leaf. In \cite[Lem. 3.3]{Kuhl} and \cite[Sec.3]{Vaq3}, it is proved that  $\mu$ is an inner node if and only if $\res(\mu)$ is an inner node. As a consequence, $\mu$ is an infinite leaf if and only if $\res(\mu)$ is an infinite leaf.

In this section, we review some basic facts on the fibers of the mapping 
$\tkb\stackrel{\res}\twoheadrightarrow\tk$. In section \ref{secKh}, we shall prove that 
$\tkh\stackrel{\res}\twoheadrightarrow\tk$ is an isomorphism of posets.

\subsection{Quasi-cuts in $\gq$ and ultrametric balls in $\kb$}\label{subsecCuts}
A \emph{quasi-cut} in a totally ordered set $X$ is a pair $D=(D^L,D^R)$ of subsets of $X$ such that 
$$
D^L\le D^R\quad\mbox{ and }\quad D^L\cup D^R=X.
$$
The condition $D^L\le D^R$ means that $\al\le \be$ for all $\al\in D^L$, $\be\in D^R$. 
Note that $D^L\cap D^R$ consists of at most one element.

Every element $\dta\in\La\infty$ determines a quasi-cut $\dd=\left(\dd^L,\dd^R\right)$ in $\gq$, defined as:
$$
\dd^L=\left(\gq\right)_{\le \dta}=\left\{\al\in\gq\mid \al\le \dta\right\},\qquad\dd^R=\left(\gq\right)_{\ge \dta}=\left\{\al\in\gq\mid \al\ge \dta\right\}
$$

Consider the equivalence relation on $\La\infty$ determined by the following conditions, where  $\dta,\dta'\in\La$: 
$$
\dta\scut\dta'\ \sii\ \dd=\ddp,\qquad \dta\not\scut\infty,\qquad \infty\scut\infty.
$$
Note that, if $\dta'<\dta$ in $\La$, then  
$$
\dta\not\scut\dta'\ \sii\  \dta'\le\al\le\dta \quad\mbox{ for some }\al\in\gq,
$$
because either $\left(\gq\right)_{\le \dta'}\subsetneq \left(\gq\right)_{\le \dta}$ or $\left(\gq\right)_{\ge \dta}\subsetneq \left(\gq\right)_{\ge \dta'}$. The next result follows.

\begin{lemma}\label{commen}
For all $\dta\in\gq\infty$, $\dta'\in\La\infty$, we have $\dta\scut\dta'$ if and only if $\dta=\dta'$. 
\end{lemma}

For all $a\in \kb$, $\dta\in\La\infty$, consider the closed ultrametric ball in $\kb$ with center $a$ and radius $\dta$:
$$
\bad=\left\{b\in\kb\mid \vb(a-b)\ge\dta\right\}.
$$

Note that $B(a,\infty)=\{a\}$ and $\bad=B(b,\dta)$ for all $b\in \bad$.

\subsection{Extensions of $v$ to $\kbx$}\label{subsecOm}
%Recall that we have fixed an extension $\vb$ of $v$ to $\kb$.

The well-specified extensions of $\vb$ to $\kbx$ are the monomial valuations $\om=\om_{a,\dta}$, for some $a\in \kb$ and $\dta\in \La\infty$. They act as follows on $(x-a)$-expansions 
$$%\begin{equation}\label{defom}
\om\left(\sum\nolimits_{0\le s}a_s(x-a)^s\right) = \min\left\{\vb(a_s)+s\dta\mid0\le s\right\}.
$$%&\end{equation}

This valuation $\om_{a,\dta}$ is commensurable over $\vb$ if and only if $\dta\in\gq\infty$. 

The proof of the following observation is an easy exercise.

\begin{lemma}\label{balls2}
Take any $a,b\in \kb$ and $\dta,\dta'\in\La\infty$. Then, 
$$\om_{b,\dta'}\le\om_{a,\dta}\ \sii\ \dta'\le\min\{\dta,\vb(a-b)\}
\ \sii\ B(b,\dta')\supset\bad \ \mbox{ and }\ \dta'\le\dta.
$$
\end{lemma}

Let us characterize equality and equivalence of well-specified valuations on $\kbx$ in terms of their associated ultrametric balls.

\begin{lemma}\label{OmBalls}
Take any $a,b\in \kb$ and $\dta,\dta'\in\La\infty$. then,
$$
\begin{array}{l}
 \om_{b,\dta'}\sim\om_{a,\dta}\ \sii\  B(b,\dta')=\bad\  \mbox{ and }\ \dta\scut\dta'.\\
\om_{b,\dta'}=\om_{a,\dta}\ \sii\ B(b,\dta')=\bad\  \mbox{ and }\ \dta=\dta'.
\end{array}
$$
\end{lemma}
 
\begin{proof}
Denote $\om=\om_{a,\dta}$, $\om'=\om_{b,\dta'}$.
The characterization of the equality  $\om'=\om$ follows from Lemma \ref{balls2}.
Let us prove the characterization of the equivalence  $\om'\sim\om$. 

If one of the valuations is commensurable, then $\om'\sim\om$ if and only if $\om'=\om$ \cite[Sec. 7.1]{VT}. By Lemma \ref{commen}, if one of the radii $\dta,\dta'$ belongs to $\gq\infty$, then $\dta\scut \dta'$ if and only if $\dta=\dta'$. Thus, the  characterization of the equivalence follows from the  characterization of the equality, in this case.

Let us suppose that both valuations are incommensurable over $\vb$; or equivalently, neither of the two radii $\dta,\dta'$ belong to $\gq\infty$. The value groups are the subgroups generated by $\gq$ and $\dta$, $\dta'$, respectively:
$$
\g_{\om}=\gen{\gq,\dta},\qquad \g_{\om'}=\gen{\gq,\dta'}.
$$
If $\om \sim \om'$, then there is an order-preserving isomorphism $\iota\colon  \g_{\om}\ism\g_{\om'}$ such that $\om'=\iota\circ\om$. Since the restriction to $\kb$ of both valuations is $\vb$, this isomorphism $\iota$ must act as the identity on $\gq$. This implies that $\dta\scut \iota(\dta)$. On the other hand, since neither $\dta'$ nor $\iota(\dta)$ belong to $\gq$, from the equality 
$$
\dta'=\om'(x-b)=(\iota\circ\om)(x-b)=\min\{\iota(\dta),\vb(a-b)\}
$$
we deduce that  $\iota(\dta)=\dta'<\vb(a-b)$.  Hence, $B(b,\dta')=B(a,\dta')=\bad$, because $\left(\gq\right)_{\ge \dta}=\left(\gq\right)_{\ge \dta'}$.

Conversely, suppose that $\vb(a-b)\ge\dta$ and $\dta\scut\dta'$. Then, the group isomorphism
$$
\iota\colon \gen{\gq,\dta}\lra\gen{\gq,\dta'},\qquad  \al+\be\dta\longmapsto \al+\be\dta'
$$
preserves the ordering and acts as the identity on $\gq$. Finally, $\om'=\iota\circ\om$, because for all $c \in \kb$ we have
\begin{align*}
 \om'(x-c)=&\min\{\dta',\vb(b-c)\}=\min\{\iota(\dta),\vb(c-b)\}\\=&\min\{\iota(\dta),\vb(a-c)\}=(\iota\circ\om)(x-c).
\end{align*}
This ends the proof of the lemma.
\end{proof}\e

Clearly, $x-a$ \emph{defines} $\om=\om_{a,\dta}$ in the terminology of Section \ref{subsecWell}. 

If $\om$ is an inner node of $\tkb$ ($\dta<\infty$), we may easily characterize the set $\kp(\om)$. 

\begin{lemma}\label{kpom}
Let $\om=\om_{a,\dta}$ for some $a\in\kb$, $\dta\in\La$, and let $b\in\kb$. The following conditions are equivalent.
\begin{enumerate}
\item[(i)] \ $x-b\in\kpo$.
\item[(ii)] \ $x-b$ is $\om$-minimal.
\item[(iii)] \ $\ino(x-b)$ is not a unit in $\ggo$. 
\item[(iv)] \ $b\in \bad$.
\item[(v)] \ $\om(x-b)=\dta$.
\item[(vi)] \ $\om=\om_{b,\dta}$.
\end{enumerate} 	
\end{lemma}	
 
\begin{proof}
Clearly, (iv), (v) and (vi) are equivalent. On the other hand, the implications (vi)$\Rightarrow$(i)$\Rightarrow$(ii)$\Rightarrow$(iii) are obvious. 

Thus, it sufficies to show that (iii) implies (iv). Indeed,
$$
b\notin \bad\imp\vb(a-b)<\dta\imp x-b\sim_\om a-b,
$$
and this implies that $\ino(x-b)=\ino(a-b)$ is a unit in $\ggo$.
\end{proof}\e

The finite leaves $\lfin(K)\subset\tk$ are the monomial valuations $\om=\om_{a,\infty}$, whose support is $(x-a)\kbx$.
The infinite leaves $\li(K)\subset\tk$ are determined by descendent chains of ultrametric balls with empty intersection \cite[Sec. 4]{Vaq3}.

\subsection{Extended value of polynomials with respect to a valuation}\label{subsecLevel}

Denote by $\N\subset \N_0$ the sets of positive and non-negative integers, respectively.
For all $b\in \N_0$, consider the linear differential operator $\partial_b$ on $\kx$, defined by Taylor's formula:
$$
f(x+y)=\sum_{0\le b}\pb f\,y^b, \quad\ \forall\,f\in\kx,
$$
where $y$ is another indeterminate. 

%Let $f\in\kx$ be a polynomial of positive degree.  Denote $$\ml(f)=\mbox{least $b\in\N$ such that }\prt{b}f\ne0. $$

%Clearly, $\ml(f)=1$ if $\chr(K)=0$. If $\chr(K)=p$, then $\ml(f)=p^r$ is the largest power of $p$ such that $f$ belongs to $K[x^{p^r}]$. \e

We need some concepts that were introduced in \cite{Dec}.\e

\noindent{\bf Definition. }Let $\mu$ be a well-specified valuation on $\kx$ and take $f\in\kx\setminus K$. 

If $\mu(f)=\infty$, we define $\epm(f)=\infty$. Otherwise, we define
	$$
	\epm(f)=\max\left\{\dfrac{\mu(f)-\mu\left(\pb{f}\right)}b\ \Big|\ b\in\N\right\}\in\La.
	$$
We say that $\epm(f)$ is the \emph{extended value} of $f$ with respect to $\mu$. 

For all $a\in K$, we agree that $\ep_\mu(a)=-\infty$. \e

If $\mu(f)<\infty$, let $I_\mu(f)$ be the set of all $b\in\N$ for which $\mu(f)-\mu(\pb{f})=b\,\epm(f)$. 

If $\mu(f)=\infty$ and $f$ is irreducible, we agree that $I_\mu(f)=\{\ml(f)\}$, where $\ml(f)$ is the least $b\in\N$ such that $\pb f\ne0$. 

Otherwise, the set  $I_\mu(f)$ is not defined.

Novacoski found an interesting interpretation of $\epm(f)$ in \cite[Prop. 3.1]{N2019}.

\begin{proposition}\label{ep=max}
Let $\mu$ be a well-specified valuation on $\kx$ and take $f\in\kx\setminus K$. Let $Z(f)$ be the set of roots of $f$ in $\kb$. 
For any extension $\om$ of $\mu$ to $\overline{K}[x]$, we have
$$
\epm(f)=\max\{\om(x-c)\mid c\in Z(f)\}.
$$
\end{proposition}

\noindent{\bf Definition. }
Let $\mu$ be a well-specified valuation on $\kx$.
A monic $Q\in\kx$ is said to be an \emph{abstract key polynomial} of $\mu$ if for all $f\in\kx$ we have
$$
\deg(f)<\deg(Q)\imp \epm(f)<\epm(Q).
$$

For an abstract key polynomial, the truncation function $\mu_Q$ is a valuation on $\kx$ such that $\mu_Q\le\mu$. We recall that $\mu_Q$ acts as follows on $Q$-expansions. 
$$
\mu_Q\left(\sum\nolimits_{0\le \ell}a_\ell Q^\ell\right)= \min\left\{\mu\left(a_\ell Q^\ell\right)\mid 0\le\ell\right\}.
$$

\noindent{\bf Definition. }Let $\mu$ be a well-specified valuation. A polynomial $f\in\kx$ is said to be \emph{$\epm$-maximal} if $\ep_\mu(f)=\max\left(\ep_\mu\left(\kx\right)\right)$.\e

The following results follows immediately from Lemma \ref{defining} and \cite[Cor. 2.22]{AFFGNR}.

\begin{proposition}\label{def=maxep}
Let $\mu$ be a well-specified valuation on $\kx$. For all monic $\phi\in \kx$, the following conditions are equivalent. 
\begin{enumerate}
	\item [(i)] $\phi$ is irreducible and it defines $\mu$.
	\item [(ii)] $\phi$ is an $\epm$-maximal polynomial of minimal degree.
\item [(iii)] $\phi$ is an abstract key polynomial of $\mu$ such that $\mu_\phi=\mu$.   
\end{enumerate}
\end{proposition}

Actually, in Lemma \ref{epphi} below we shall characterize all $\epm$-maximal polynomials.

Finally, let us quote a trivial property of the extended value. 

\begin{lemma}\label{e<e}
 Let $\nu<\mu$ be two well-specified valuations, and let $\ty(\nu,\mu)=[\varphi]_\nu$ be the corresponding tangent direction. Then, $\epn(\varphi)<\epm(\varphi)$. 
\end{lemma}

This follows immediately from
 $$
\nu(\varphi)<\mu(\varphi),\qquad \nu(\pb \varphi)=\mu(\pb \varphi),\ \forall\,b\in\N.
$$

\subsection{Up and down}\label{subsecUD}

For any ultrametric ball $B\subset\kb$, we define its $K$-degree as
$$
\deg_KB=\min\{\deg_Ka\mid a\in B\}.
$$
A pair $(a,\dta)\in\kb\times\La$ is said to be \emph{$K$-minimal} if $\deg_KB(a,\dta)=\deg_Ka$.   

The following result on the restriction mapping $\tkb\stackrel{\res}{\to}\tk$ is a consequence of the results in \cite{APZ} and \cite{PP}. Another proof can be found in \cite[Thm. 1.1]{N2019}.

\begin{proposition}\label{minpair}
	Let $(a,\dta)\in\kb\times\La$ be a $K$-minimal pair and $\mu=\res\left(\om_{a,\dta}\right)$. 
	Then, $\irr_K(a)$ is a key polynomial for $\mu$, of minimal degree and 
	$$
	\mu(f)=\vb(f(a))\quad \mbox{for all }\ f\in \kx \ \mbox{ with }\ \deg(f)<\deg_Ka.
	$$	
	
\end{proposition}

\begin{corollary}\label{degk}
	For all $(a,\dta)\in\kb\times\La$, we have
	$\deg_K \bad=\deg(\res(\om_{a,\dta}))$.
\end{corollary}

For a given well-specified valuation on $\kx$, our aim is to understand how to construct the extensions of $\mu$ to $\kbx$. 
To this end, the following result plays an essential role.

\begin{proposition}\cite[Prop. 3.3]{Vaq3}\label{unitunit}
	For an inner node $\om$ of $\tkb$, let $\mu=\res(\om)$. 
	For all $f\in\kx$, $\inn_\mu f$ is a unit in $\ggm$ if and only if   \,$\inn_{\om} f$ is a unit in $\gg_{\om}$.
\end{proposition}

All ultrametric balls in $\kb$ inducing $\mu$ have the same radius. This is proven in \cite[Sec.4]{MMS} and \cite[Sec. 3]{Vaq3}. Let us provide a short proof.  

\begin{lemma}\label{radiusup}
	Let $\mu$ be an inner node of $\tk$. Take any $\phi\in\kpm$ and denote $\dta=\epm(\phi)$. If $\om\in\tkb$ satisfies $\res(\om)=\mu$, then  $\om=\om_{a,\dta}$ for some $a\in Z(\phi)$.
\end{lemma}

\begin{proof}
	Since $\inn_\mu \phi$ is not a unit in $\ggm$, Proposition \ref{unitunit} shows that $\inn_{\om}\phi$ is not a unit
	in $\gg_{\om}$. Thus, there exists $a\in Z(\phi)$ such that $\inn_{\om}(x-a)$ is not a unit. By Lemma \ref{kpom}, $\om=\om_{a,\ga}$ for $\ga=\om(x-a)$.
	Since $\om(x-b)\le\ga$ for all $b\in\kb$, Proposition \ref{ep=max} shows that $\ga=\dta$.
\end{proof}\e

Actually, this common radius is the maximal value of $\epm$ in $\kx$. 

\begin{lemma}\label{epphi}
	Let $\mu$ be an inner node of $\tk$. A polynomial $f\in\kx$ is $\epm$-maximal if and only if $\,\inn_\mu f$ is not a unit in $\ggm$. 
\end{lemma}

\begin{proof}
	Take $\phi\in\kpm$ of minimal degree in this set. By Lemma \ref{defining}, $\phi$ defines $\mu$. By Proposition \ref{def=maxep}, $\phi$ is an abstract key polynomial of $\mu$ and  $\epm(\phi)=\max\left(\epm\left(\kx\right)\right)$, so that $\epm(f)\le\epm(\phi)$.
	If $f$ is not $\epm$-maximal, then $\epm(f)<\epm(\phi)$ and $\,\inn_\mu f$ is a unit in $\ggm$ by \cite[Cor. 2.8]{AFFGNR}.
	
	Conversely, suppose that $\,\inn_\mu f$ is a unit in $\ggm$. For $\dta=\epm(\phi)$, take any $a \in Z(\phi)$ such that $\om=\om_{a,\dta}$ restricts to $\mu$. By Proposition \ref{unitunit}, $\inn_{\om} f$ is a unit in $\gg_{\om}$, so that
	$\inn_{\om} (x-c)$ is a unit in $\gg_{\om}$ for all $c\in Z(f)$. By Lemma \ref{kpom}, $\om(x-c)<\dta$, so that $\epm(f)<\dta$ by Proposition \ref{ep=max}.
\end{proof}%\e

\begin{corollary}\label{twocases}
	Let $\mu\in \tinnk$ be an inner node.
	\begin{enumerate}
		\item[(i)] All key polynomials for $\mu$ are $\epm$-maximal. 
		\item[(ii)] Let $\nu\in\tk$ such that $\mu<\nu$. Then, all polynomials $f\in\kx$ such that $\mu(f)<\nu(f)$ are $\epm$-maximal. 
	\end{enumerate}
\end{corollary}

\begin{proof}
	If $\phi\in\kpm$, then $\inn_\mu \phi$ is an homogeneous prime element in $\ggm$. Thus, it is not a unit. By Lemma \ref{epphi}, $\phi$ is $\epm$-maximal.
	
	Let $\phi\in\kpm$ be such that $\ty(\mu,\nu)=[\phi]_{\mu}$. The condition $\mu(f)<\nu(f)$ is equivalent to $\phi\mmu f$. Thus, $\inn_\mu f$ cannot be a unit. By Lemma \ref{epphi}, $f$ is $\epm$-maximal.
\end{proof}

\section{Lifting valuations from $\kx$ to $\khx$}\label{secKh}

\subsection{Finite leaves}\label{subsecLfin}
For any given field $L$ and a monic irreducible polynomial $f\in\irr(L)$, we denote by $L_f$ the simple finite field extension $L[x]/(f)$. 

The set $\lfin(K)$ of finite leaves of $\tk$ may be parametrized as
$$
\lfin(K)=\left\{(f,\nb)\mid f\in\irr(K),\ \nb\ \mbox{ valuation on $K_f$ extending }v \right\},
$$
where we identify each pair $(f,\nb)$ with the following valuation with support $f\kx$:
$$
\nu\colon \kx\longtwoheadrightarrow K_f\stackrel{\nb}\lra\gq\infty.
$$

In the henselian case, the valuation $\vh$ has a unique extension to $\kb$ and we get a bijective mapping 
$$
\irr(\kh)\lra\lfin(\kh),\qquad F\longmapsto v_F,
$$
where the valuation $v_F$ on $\khx$ acts as follows
$$
v_F\colon \khx\lra \gq\infty,\qquad g\longmapsto v_F(g)=\vb\left(g(a)\right),
$$
for any root $a\in \kb$ of $F$. By the henselian property, this construction does not depend on the choice of $a$ in the set $Z(F)\subset \kb$ of all roots of $F$. 

For each polynomial $f\in\irr(K)$, let us describe the extensions of $v$ to the simple extension $K_f$. 
Since $K^h/K$ is a separable extension, the factorization of $f$ into a product of monic irreducible polynomials in $K^h[x]$ takes the form
$$
f=F_1\cdots F_r,\qquad F_1,\dots,F_r\in\irr(K^h),
$$
pairwise different. Let $w_{F_i}$ be the restriction to $\kx$ of the valuation $v_{F_i}$ and let $\bar{w}_{F_i}$ be the corresponding valuation on $K_f$. The following result is classical \cite[Sec. 17]{endler}, \cite[Sec. 3]{VT}.

\begin{theorem}\label{endler}
	There are $r$ extensions  of $v$ to $K_f$, given by $\bar{w}_{F_1},\dots,\bar{w}_{F_r}$. 
\end{theorem}

In particular, there are $r$ finite leaves of $\tk$ with support $f\kx$, given by
$$
w_{F_1}=(f,\bar{w}_{F_1}),\dots, w_{F_r}=(f,\bar{w}_{F_r}).
$$

\begin{corollary}\label{lfin}
The following mappings are both bijective $$\irr(\kh)\lra\lfin(\kh)\stackrel{\res}\lra\lfin(K).$$
\end{corollary}

Indeed, the first mapping is bijective, and Theorem \ref{endler} shows that the composition 
$\irr(\kh)\to\lfin(K)$, mapping $F$ to $w_F$, is bijective too.

We may use this bijection to associate an irreducible polynomial over $\khx$ to any key polynomial for a valuation on $\kx$.\e

\noindent{\bf Definition. }Let $\phi\in\kx$ be a key  polynomial for a valuation $\mu$ on $\kx$. The  \emph{distinguished factor} of the pair $(\mu,\phi)$ is the irreducible factor $\fmp\in\irr(\kh)$ of $\phi$  associated to the finite leaf $[\mu;\,\phi,\infty]\in\lfin(K)$ through the above bijection. \e

In other words, $[\mu;\,\phi,\infty]=w_{\fmp}$.
By Theorem \ref{endler}, $\fmp$ is characterized as the unique  $F\in\irr(\kh)$ satisfying the following property:
$$
g\in\kx,\quad \deg(g)<\deg(\phi) \imp \mu(g)=\vb(g(a)),
$$
for some (any) root $a\in Z(F)$.

\subsection{Rigidity}
Let $\mu\in\tinnk$ be an inner node and take $\phi\in\kpm$. Let $\dta=\epm(\phi)=\max\left(\epm(\kx)\right)$.
We have seen in section \ref{subsecUD} that for some roots $a\in Z(\phi)$, we have $\res(\om_{a,\dta})=\mu$.

However, not all roots of $\phi$ have this property.  Which are the ``good" roots? 

The following result answers this question when $\phi$ is a key polynomial for $\mu$ of minimal degree. 
%For a key polynomial $\phi$ of minimal degree, a criterion is found  in \cite[Lem. 6.1]{MMS}.  

\begin{proposition}\label{up0}
	Let $\mu$ be an inner node of $\tk$. Take $\phi\in\kpm$ defining $\mu$ and denote $\dta=\epm(\phi)$. Then, for all $a\in Z(\phi)$, we have
$$
\res(\om_{a,\dta})=\mu\ \sii\ a\in Z(\fmp).
$$
In this case, $(a,\dta)$ is a $K$-minimal pair. 
\end{proposition}

\begin{proof}
Let us denote $F=\fmp$. Recall that $w_F=[\mu;\,\phi,\infty]$.

	Suppose that $\res(\om_{a,\dta})=\mu$. By Corollary \ref{degk}, $$\deg_K\bad=\deg(\mu)=\deg(\phi)=\deg_Ka.$$ Thus, $(a,\dta)$ is a $K$-minimal pair. By Proposition \ref{minpair}, 
$$
\mu(f)=\vb(f(a))\quad \mbox{for all }\ f\in \kx \ \mbox{ with }\ \deg(f)<\deg_Ka.
$$	
By Theorem \ref{endler}, this implies that $a$ is a root of $F$. 

Conversely, suppose that $a\in Z(F)$. By Lemma \ref{radiusup}, there exists a root $b\in Z(\phi)$ such that $\res(\om_{b,\dta})=\mu$. By the implication we have just proven, $b$ is necessarily a root of $F$. Hence, there exists $\sg\in\aut(\kb/\kh)$ such that $\sg(a)=b$.
Since $\vb\circ \sg=\vb$, we deduce that $\om_{a,\dta}\circ \sg= \om_{b,\dta}$. Hence,
\begin{equation}\label{equalres}
f\in\khx\ \imp\ \sg(f)=f\ \imp\ \om_{a,\dta}(f)= \om_{b,\dta}(f).
\end{equation}
In particular, $\res(\om_{a,\dta})=\res(\om_{b,\dta})=\mu$.
\end{proof}\e

Proposition \ref{up0} may be considered a reinterpretation of \cite[Lem. 6.1]{MMS}, in the light of Theorem \ref{endler}. 

\begin{corollary}\label{twsBij}
The restriction mapping $\twskh\to\twsk$ is bijective.
\end{corollary}

\begin{proof}
If $\mu\in \twsk$ is a finite leaf, it has a unique preimage in $\twskh$ by Corollary \ref{lfin}.

If $\mu\in \twsk$ is an inner node, then Proposition \ref{up0} shows that all extensions of $\mu$ to $\kbx$ have the same restriction to $\khx$, as shown in (\ref{equalres}). Thus, $\mu$ has a unique preimage in $\twskh$ too.
\end{proof}  \e

For all $\mu\in\twsk$, let us denote by $\muh$ the unique preimage of $\mu$ in $\tkh$.

\begin{theorem}\label{twsIso}
The restriction mapping $\twskh\to\twsk$ is an isomorphism of posets.
\end{theorem}

Obviously, the restriction mapping preserves the ordering; thus, we need only to show that its inverse mapping
$$
\res^{-1}\colon \twsk\lra\twskh,\qquad \mu\longmapsto\muh 
$$
preserves the ordering too. 

By the main theorem of MacLane-Vaqui\'e, if $\nu,\mu\in\twsk$ satisfy $\nu<\mu$, then $\mu$ may be constructed by a finite number of augmentation steps (ordinary or limit) starting with $\nu$. Therefore, Theorem \ref{twsIso} is a consequence of the following lemma. 

\begin{lemma}\label{AugmIso}
Let $\nu\stackrel{\phi,\ga}\lra\mu$ be an augmentation step  in $\twsk$, either ordinary or limit. Then, $\nuh<\muh$ in $\twskh$. 
\end{lemma}

\begin{proof}
Let $\dta'=\epn(\phi)$ and $\dta=\epm(\phi)$.

If $\mu$ is a finite leaf, then $\ga=\dta=\infty$ and $\mu=w_F$ for some irreducible factor $F$ of $\phi$ in $\khx$. Clearly, the extensions of $\mu$ to $\tkb$ are the finite leaves of the form $\om_{a,\infty}$ for $a\in Z(F)$. 

If $\mu$ is an inner node, then $ \phi$ is a key polynomial for $\mu$ of minimal degree, by Lemma \ref{defining}. By Proposition \ref{up0}, there exists an irreducible factor $F$ of $\phi$ in $\khx$ such that the extensions of $\mu$ to $\tkb$ are of the form $\om_{a,\dta}$ for $a\in Z(F)$. 

In all cases, we may fix a root $a\in Z(\phi)$ such that  $\mu=\res(\om_{a,\dta})$. Denote 
$$\om'=\om_{a,\dta'}, \qquad \om=\om_{a,\dta},\qquad \nu'=\res(\om')\in\tk.$$ 
By Proposition \ref{ep=max}, $\ep_{\nu'}(\phi)=\dta'=\epn(\phi)$.\e

\noindent{\bf Claim. }$\dta'<\dta$ and $\nu'=\nu$.\e
  
From this claim, the lemma follows immediately. Indeed, 
since $\dta'<\dta$, we have $\om'\le\om$, and hence, $\nu'\le\mu$. Therefore,
$$
\nuh=(\nu')^h=\res_{\kb/\kh}(\om')\le\res_{\kb/\kh}(\om)=\muh.
$$
Also, since $\nu\ne\mu$, we have $\nuh\ne \muh$ by Corollary \ref{twsBij}.  \e
  
Let us first prove the Claim in the ordinary augmentation case. Suppose that $\mu=[\nu;\,\phi,\ga]$  for some $\phi\in\kpn$, $\ga\in\La\infty$, $\ga>\nu(\phi)$. 
  Since $\ty(\nu,\mu)=[\phi]_\nu$, Lemma \ref{e<e} shows that $\dta'<\dta$. In particular, $\nu'<\mu$.

Since $\tk$ is a tree, the set of valuations bounded above by $\mu$ is totally ordered. Hence, we have two possibilities:
$$
\nu'\le\nu<\mu,\quad\mbox{ or }\quad \nu\le\nu'<\mu.
$$

Suppose that $\nu'<\nu<\mu$. Let $\ty(\nu',\nu)=[\varphi]_{\nu'}$. Since $\varphi\in\kp(\nu')$ and $\nu(\phi)<\mu(\phi)$, Corollary \ref{twocases}  shows that $\varphi$ is $\ep_{\nu'}$-maximal and $\phi$ is $\epn$-maximal.
By applying Lemma \ref{e<e} we get a contradiction:
$$
\dta'=\ep_{\nu'}(\phi)\le\ep_{\nu'}(\varphi)<\epn(\varphi)\le\epn(\phi)=\dta'.
$$

Suppose that $\nu<\nu'<\mu$. Then, $\ty(\nu,\nu')=\ty(\nu,\mu)=[\phi]_{\nu}$ and Lemma \ref{e<e} leads  to a contradiction too:
$$
\dta'=\epn(\phi)<\ep_{\nu'}(\phi)=\dta'.
$$
Therefore, $\nu=\nu'$.
This ends the proof of the lemma for ordinary augmentations.\e

Let us prove the Claim in the limit augmentation case. Suppose that $\mu=[\cc;\,\phi,\ga]$ for some essential continuous family  $\cc=\left(\ri\right)_{i\in A}\subset \tinnk$ whose initial valuation is $\nu$, some limit key polynomial $\phi\in\kpi(\cc)$ and some value $\ga=\mu(\phi)\in\La\infty$, $\ga>\rho_i(\phi)$ for all $i\in A$. 

Since $\deg\left(\pb \phi\right)<\deg(\phi)$, these polynomials are $\cc$-stable for all $b\in\N$. Since $\pb \phi=0$ for $b>\deg(\phi)$, there exists $i_0\in A$ such that
$$
\ri(\pb \phi)=\mu(\pb \phi)\quad\mbox{ for all }b\in\N,\quad\mbox{ for all }i\ge i_0.
$$ 

Since $\deg(\nu)=\deg(\rho_{i_0})$, the step $\nu\to\rho_{i_0}$ is covered by a single ordinary augmentation \cite[Lem 4.5]{VT}. As we proved above, this implies $\nuh<\rho_{i_0}^h$. Thus, it suffices to show that $\rho_{i_0}^h<\muh$, so that we may assume that $\nu=\rho_{i_0}$. 

Since $\ri(\phi)<\mu(\phi)$ for all $i\in A$, we deduce that $\dta'=\epn(\phi)<\epm(\phi)=\dta$.

As above, this implies that
$$
\nu'\le\nu<\mu,\quad\mbox{ or }\quad \nu\le\nu'<\mu.
$$

If $\nu'<\nu<\mu$, we get a contradiction by exactly the same arguments as above.

Suppose that $\nu<\nu'<\mu$. Then, $\ty(\nu,\nu')=\ty(\nu,\mu)=[\chi]_{\nu}$ for some $\chi\in\kpn$. Since $\nu(\phi)<\mu(\phi)$, we have $\chi\mid_\nu\phi$, so that $\nu(\phi)<\nu'(\phi)$ as well. Since $\nu(\pb \phi)=\nu'(\pb \phi)$ for all $b\in\N$, we get a contradiction too:   
$$
\dta'=\epn(\phi)<\ep_{\nu'}(\phi)=\dta'.
$$
Therefore, $\nu=\nu'$.
This ends the proof of the lemma.
\end{proof}

\subsection{Lifting of infinite leaves}

Let us briefly review some properties of infinite leaves, mainly extracted from \cite[Sec. 4]{VT}.
Consider totally ordered families 
$$
\aa=\cfa,\quad \ri\in\tk,
$$
parameterized by a totally ordered set $A$ of indices such that the mapping $A\to\aa$, given by $i\mapsto \ri$,
is an isomorphism of totally ordered sets.
Moreover, we shall always assume that $\aa$ contains no maximal element. Therefore, all valuations $\ri$ are inner nodes of  $\tk$. 

The degree function $\deg\colon \aa\to\N$
is order-preserving \cite[Lem. 2.2]{VT}. Hence, these families fall into two radically different cases:\e

(a) \ The set $\deg(\aa)$ is unbounded in $\N$. We say that $\aa$ has \emph{unbounded degree}.\e

(b) \ There exists $i_0\in A$ such that $\deg(\ri)=\deg(\rho_{i_0})$ for all $i\ge i_0$. 
We say that $\aa$ is a \emph{continuous family} of stable degree $m(\aa)=\deg(\rho_{i_0})$.\e

We say that a nonzero $f\in\kx$ is \emph{$\aa$-stable} if for some index $i\in A$, it satisfies 
$$
\ri(f)=\rj(f), \ \mbox{ for all }\,j>i.
$$
In this case, we denote this stable value by  $\rha(f)$.  

We say that $\aa$ has a \emph{stable limit} if all polynomials in $\kx$ are $\aa$-stable. 
In this case, the stability function 
$$\rha\colon \kx\lra \La\infty$$ 
is a valuation. We say that $\rha$ is the stable limit of $\aa$ and we write $\rha=\lim\aa$.

These limit valuations  have trivial support and satisfy $\kp(\rha)=\emptyset$ \cite[Prop. 3.1]{MLV}. Thus, they are infinite leaves of the tree $\tk$. By the main theorem of MacLane-Vaqui\'e, all infinite leaves of $\tk$ arise in this way \cite[Thm. 4.3]{MLV}.

All totally ordered families of unbounded degree have a stable limit \cite[Cor. 4.4]{VT}. %, and the function $$\sval\colon \aa\lra\La,\qquad \ri\longmapsto \sval(\ri) $$is unbounded too.
However, the continuous families may have $\aa$-unstable polynomials, leading to limit augmentations.

The following lemma shows that stable limits behave well with respect to restriction. The proof follows immediately from the definitions.

\begin{lemma}\label{LimRes}
Let $\aa$ be a totally ordered family in $\tkh$ admitting stable limit. Then, $\res(\aa)$ is a totally ordered family in $\tk$ admitting stable limit and
$$
\lim\res(\aa)=\res(\lim\aa).
$$ 
\end{lemma}

%\e

\noindent{\bf Definition. }Let $\aa=\cfa$, $\bb=\cfb$ be two totally ordered families having a stable limit.
We say that $\aa$ and $\bb$ are \emph{equivalent}, and we write $\aa\sim\bb$, if they are cofinal in each other.\e

%Since the degree mapping preserves the ordering, two equivalent families either both have unbounded degree, or both have the same stable degree.

%The set $\li(K)$ of infinite leaves splits into two disjoint subsets, gathering the limits of unbounded, or bounded, degree families.

 The set of infinite leaves is parameterized by the equivalence classes of totally ordered families of valuations having a stable limit.
 
\begin{proposition}\cite[Prop. 4.9]{VT}\label{equiv=lim}
Let $\aa$, $\bb$ be two totally ordered families having a stable limit. Then,
$$
\aa\sim\bb\ \sii\ \lim\aa=\lim\bb.
$$
\end{proposition}

As a consequence, Theorem A follows immediately from Theorem \ref{twsIso}.\e

\noindent{\bf Theorem A. }{\it The restriction mapping $\tkh\to\tk$ is an isomorphism of posets.}\e

\noindent{\bf \quad Proof. }Let us first show that the restriction mapping $\li(\kh)\twoheadrightarrow\li(K)$ is bijective. Suppose that $\eta=\lim\aa$, $\xi=\lim\bb$ are two infinite leaves in $\li(\kh)$ with the same restriction to $\li(K)$.
By Theorem \ref{twsIso}, we have
$$
\aa\sim\bb\ \sii\ \res(\aa)\sim\res(\bb).
$$
If we combine this with Lemma \ref{LimRes}, we obtain \begin{align*}
\res(\eta)=\res(\xi)&\ \imp\ \lim \res(\aa)=\lim \res(\bb)\ \imp\ \res(\aa)\sim\res(\bb)\\&\ \imp\ \aa\sim\bb\ \imp\ \eta=\xi.
\end{align*}
Therefore, the restriction mapping $\li(\kh)\to\li(K)$ is bijective. By Corollary \ref{twsBij}, we deduce that the restriction mapping $\tkh\to\tk$ is bijective too.

Finally, let us see that the inverse mapping
$$
\res^{-1}\colon \tk\lra\tkh,\qquad \mu\longmapsto \muh 
$$
preserves the ordering. 

Suppose that $\nu<\mu$ in $\tk$. If $\mu\not\in\li(K)$, then $\nu\not\in\li(K)$ and hence $\nuh<\muh$ by Theorem \ref{twsIso}. If $\mu=\lim\aa$ is an infinite leaf, there exists $\rho\in \aa$ such that $\nu<\rho<\mu$. Since $\nuh<\rho^h$, we may replace $\nu$ by $\rho$ and assume that $\nu\in \aa$. Now, consider the family $\aa^h$ formed by all liftings of the valuations in $\aa$ to $\khx$. Clearly, $\aa^h$ has a stable limit and 
$$\nuh<\eta=\lim\aa^h\in\li(\kh).$$ By Lemma \ref{LimRes}, $\res(\eta)=\mu$. Since the restriction mapping $\li(\kh)\to\li(K)$ is bijective, we have necessarily $\eta=\muh$.\hfill{$\Box$}\e

\section{Equivalence of valuations}\label{secEquiv}

If $\mu,\nu$ are two equivalent valuations on $\kx$, we write $\mu\sim\nu$.
Clearly, these valuations have isomorphic graded algebras.\e

\noindent{\bf Definition. }
Let $\tv$ be the set of equivalence classes of valuations on $\kx$ whose restriction to $K$ is equivalent to $v$. 
Let $\tvh$ be the set of equivalence classes of valuations on $\khx$ whose restriction to $\kh$ is equivalent to $\vh$. \e

We emphasize that we are not fixing any ordered abelian group containing the value groups of these valuations.

The set $\tv$ admits a natural decomposition into a disjoint union of the subsets:
$$
\tv=\tvinn\,\sqcup\,\lfin(\tv)\,\sqcup\,\li(\tv).
$$
Although we are not yet considering an ordering on $\tv$, we say that these are the subsets of \emph{inner nodes}, \emph{finite leaves} and \emph{infinite leaves}, respectively. These subsets are defined in the usual way:

\begin{itemize}
\item \ $\tvinn=\left\{[\mu]\mid \kpm\ne\emptyset\right\}$.
\item \ $\lfin(\tv)=\left\{[\mu]\mid \supp(\mu)\ne0\right\}$.
\item \ $\li(\tv)=\left\{[\mu]\mid \kpm=\emptyset,\ \supp(\mu)=0\right\}$.
\end{itemize}

It is well known how to describe the equivalence classes of valuations $\mu$ which are commensurable over $v$. The following result is proved in \cite[Sec. 7.1]{VT}. 

\begin{lemma}\label{comm}
Let $\tv^\com\subset\tv$ be the subset of equivalence classes of valuations which are commensurable over $v$. Then, the mapping
$$
\ttt\left((K,v),\gq\right)\lra\tv,\qquad \mu\longmapsto [\mu]
$$
is injective, and its image is $\tv^\com$. 
In particular, it induces bijections $$\lfin\left((K,v),\gq\right)\ism\lfin(\tv),\qquad \li\left((K,v),\gq\right)\ism\li(\tv).
$$  	
\end{lemma}

The natural restriction mapping
$$
\res\colon \tvh\longtwoheadrightarrow \tv,\qquad [\eta]\longmapsto [\res(\eta)]
$$
is well defined and onto. Indeed, suppose that $\eta$, $\xi$ are valuations on $\khx$ such that $\eta\sim\xi$. Let $\iota\colon \g_{\eta}\ism\g_{\xi}$ be an order preserving isomorphism such that $\xi=\iota\circ\eta$. Then, it is obvious that this isomorphism restricts to an isomorphism between the subgroups $\g_{\res(\eta)}\subset\g_\eta$ and  $\g_{\res(\xi)}\subset\g_\xi$, still satisfying $\res(\xi)=\iota\circ\res(\eta)$.

Our first aim in this section is to show that this mappping is a bijection. To this end, the following result is crucial.

\begin{proposition}\label{characEquiv}
Let $\g\hk\La$ be an embedding of $\g$ into a divisible ordered group. For all $\nu,\mu\in\tinn\left((K,v),\La\right)$, we have $\nu\sim\mu$ if and only if there exists $\phi \in\kpm\cap\kpn$ having minimal degree in both sets, satisfying	
$$F_{\mu,\phi}=F_{\nu,\phi}\quad\mbox{and }\quad\mu(\phi)\scut\nu(\phi).
$$
\end{proposition}

This proposition is a slight variation of \cite[Prop. 6.3]{VT}. Although the contexts are slightly different, the arguments of the proof can be mimicked.

\begin{theorem}\label{tvBij}
The restriction mapping $\tvh\to\tv$ is bijective.
\end{theorem}

\begin{proof}
The restriction mapping splits into three ``disjoint" mappings
$$
\tvh^{\op{inn}}\lra \tvinn,\qquad \lfin(\tvh)\lra\lfin(\tv),\qquad \li(\tvh)\lra\li(\tv).
$$ 
The last two mappings are bijective by Lemma \ref{comm} and Corollary \ref{twsBij}. Thus, we need only to show that the restriction mapping is injective on inner nodes.

Let $\eta$, $\xi$ be two valuations on $\khx$ admitting key polynomials, and denote $\nu=\res(\eta)$, $\mu=\res(\xi)$. We can find equivalent valuations $\eta'\sim\eta$, $\xi'\sim\xi$ such that $\eta'$ and $\xi'$ take values in a common divisible extension $\g\hk\La$ of the value group of $v$. Hence, we may assume that all four valuations $\eta$, $\xi$, $\nu$, $\mu$ take values in $\La$.  By Corollary \ref{twsBij}, $\eta=\nuh$ and $\xi=\muh$.

Suppose that $\nu\sim\mu$, and let us see that this implies $\nuh\sim\muh$. By Proposition \ref{characEquiv}, the valuations $\nu$ and $\mu$ admit  a common key polynomial $\phi$ of minimal degree, such that $F_{\nu,\phi}=F_{\mu,\phi}$ and $\nu(\phi)\scut\mu(\phi)$.

Denote by $F=F_{\nu,\phi}=F_{\mu,\phi}\in\irr(\kh)$ this irreducible factor of $\phi$. Since
$$
[\nu;\,\phi,\infty]=w_F=[\mu;\,\phi,\infty],
$$ 
we have $\nu<w_F$ and $\mu<w_F$. Since the valuations bounded above by $w_F$ are totally ordered, we have (say) $\nu\le\mu<w_F$. For all $f\in\kx$ with $\deg(f)<\deg(\phi)$ we have 
$$
\nu(f)=w_F(f)=\mu(f)\in\gq.
$$ 
If $\nu(\phi)=\mu(\phi)$, then $\nu=\mu$ because both valuations coincide on $\phi$-expansions. In this case, $\eta=\xi$ by Corollary \ref{twsBij}, and our theorem is proven.

Suppose $\ga':=\nu(\phi)<\ga:=\mu(\phi)$. Since $\ga'\scut\ga$, Lemma \ref{commen} shows that $\ga,\ga'\not\in\gq$. Also, there is no $\al\in\gq$ such that $\ga'\le\al\le\ga$. \e

\noindent{\bf Claim. }The values $\dta'=\epn(\phi)$, $\dta=\epm(\phi)$, satisfy $\dta\scut\dta'$.\e

From this claim, the theorem follows immediately. Indeed, take any  root $a\in Z(F)$ and consider $\om=\om_{a,\dta}$, $\om'=\om_{a,\dta'}$.
By Proposition \ref{up0},
$$
\res_{\kb/K}(\om')=\nu,\qquad \res_{\kb/K}(\om)=\mu.
$$ 
Finally,  $\om\sim\om'$ by Lemma \ref{OmBalls}, so that 
$$
\nuh=\res_{\kb/\kh}(\om')\sim\res_{\kb/\kh}(\om)=\muh.
$$

Let us prove the Claim. Let us denote
$$
\al_b=\nu(\pb \phi)=\mu(\pb\phi),\quad\mbox{for all }b\in\N.
$$   
Let us first show that $I_\nu(\phi)=I_\mu(\phi)$ (cf. section \ref{subsecLevel}). Take any $b\in I_\mu(\phi)$, or equivalently
$$
(c-b)\ga> c\al_b-b\al_c,\quad\mbox{ for all }c\in\N,\ c\ne b.
$$ 
The strict inequality is a consequence of the fact that $\ga\not\in\gq$, but $\al_b,\al_c\in\gq$. 
Now, exactly the same inequality must be satisfied by $\ga'$:
$$
(c-b)\ga'> c\al_b-b\al_c,\quad\mbox{ for all }c\in\N,\ c\ne b.
$$ 
Indeed, if we had $(c-b)\ga'< c\al_b-b\al_c$ for some $c\in\N$, $c\ne b$, this would imply
$$
\ga'<\dfrac{c\al_b-b\al_c}{c-b}<\ga,\quad\mbox{ or }\quad\ga<\dfrac{c\al_b-b\al_c}{c-b}<\ga',
$$
depending on $c>b$ or $c<b$. This contradicts the fact that $\ga\scut \ga'$. Thus, $b\in I_\nu(\phi)$. By the symmetry
of the argument, we deduce that $I_\nu(\phi)=I_\mu(\phi)$.\footnote{Moreover, we deduce that $\#I_\nu(\phi)=1$, but we do not use this fact in the proof.}

Therefore, for some $b\in\N$ we have 
$$
\dta=\left(\ga-\al_b\right)/b,\qquad \dta'=\left(\ga'-\al_b\right)/b.
$$
This shows that $\dta\scut\dta'$, because any inequality $\dta'<\be<\dta$, with $\be\in \gq$, leads to $\ga'<b\be+\al_b<\ga$, with $b\be+\al_b\in\gq$, which is impossible.
\end{proof}

\subsection{Tree structure on $\tv$} 
It is well-known that we can embed $\g$ into some universal ordered abelian group $\La$ such that every valuation on  $\kx$ is equivalent to some valuation taking values in $\La$. For instance, in \cite{csme} an explicit ordered real vector space extending $\g$ is constructed with such a universal property. 

Therefore, if $\tk=\ttt((K,v),\La)$, then the quotient set $\tk/\!\!\sim$ serves as a model of $\tv$.
Lemma \ref{inherite} below shows that this quotient  inherits the tree structure of  $\tk$. %More precisely, the following property holds:
In order to prove this result, we need to show that an equivalence class is totally ordered and  closed under the formation of intervals. That is, two equivalent valuations have the following property.

\begin{lemma}\label{EquivInterval}
Let $\nu,\mu\in \tk$ such that $\nu\sim\mu$. Then $\nu$ and $\mu$ are comparable and the interval $[\nu,\mu]$ (if $\nu\le\mu$) is contained in their equivalence class $[\nu]=[\mu]$.	
\end{lemma}

\begin{proof}
If $\nu$ or $\mu$ are commensurable over $v$, then $\nu=\mu$ and $[\mu]=\{\mu\}$. The statement of the lemma is obvious in this case.

Suppose that $\nu$ and $\mu$ are incommensurable and $\nu\ne\mu$. By Proposition \ref{characEquiv}, there exists a common key polynomial of minimal degree $\phi\in\kpn\cap \kpm$ such that $F_{\nu,\phi}=F_{\mu,\phi}$ and $\nu(\phi)\scut\mu(\phi)$. As we saw along the proof of Theorem \ref{tvBij}, this implies that $\nu$ and $\mu$ are comparable, because both are bounded above by $w_F$.

Suppose that $\nu<\mu$. For all $f\in\kx$ with $\deg(f)<\deg(\phi)$ we have $\nu(f)=w_F(f)=\mu(f)$. Therefore,
 $\mu=[\nu;\,\phi,\mu(\phi)]$, because the two valuations coincide on $\phi$-expansions.

Now, consider any valuation $\nu<\eta<\mu$. By \cite[Lem. 2.7]{MLV}, $\eta=[\nu;\,\phi,\nu(\phi)]$ and $\nu(\phi)<\eta(\phi)<\mu(\phi)$. Thus, $\phi$ is a key polynomial for $\eta$ of minimal degree such that $F_{\eta,\phi}=F_{\nu,\phi}$ and $\nu(\phi)\scut\eta(\phi)$. By Proposition \ref{characEquiv}, $\eta\sim\nu$.  
\end{proof}\e

\begin{lemma}\label{inherite}
Let $\nu,\mu\in \tk$ such that $[\nu]\ne[\mu]$ and $\nu<\mu$. Then, for all $\nu'\in[\nu],\ \mu'\in[\mu]$, we have $\nu'<\mu'$.
\end{lemma}

\begin{proof}
	Let us first show that $\nu<\mu'$ for all $\mu'\in[\mu]$. 
By Lemma \ref{EquivInterval}, $\mu$ and $\mu'$ are comparable and the interval $[\mu',\mu]$ is contained in the equivalence class $[\mu]$. If $\mu\le \mu'$ we have automatically $\nu<\mu\le\mu'$. If $\mu'\le\mu$, then $\nu$ and $\mu'$ are comparable because they have $\mu$ as a common upper bound.  By Lemma \ref{EquivInterval}, $\mu'\le\nu$ would imply $\mu'\sim\nu$, against our assumption. Thus, necessarily $\nu<\mu'$.  
	
Now, a completely symmetrical argument applied to the inequality $\nu<\mu'$ shows that $\nu'<\mu'$ for all $\nu'\sim\nu$.
\end{proof}\e

Consider a universal ordered abelian group $\La$ such that $\tkh=\ttt\left((\kh,\vh),\La\right)$ contains all equivalence classes of valuations on $\khx$ restricting to $\vh$. Obviously, $\tk$ contains as well  all equivalence classes of valuations on $\kx$ restricting to $v$. In this way, we get a tree structure on $\tv$ and $\tvh$, by identifying them with the quotient sets $\tk/\!\!\sim$, $\tkh/\!\!\sim$, respectively.\e

\noindent{\bf Theorem B. }{\it The restriction mapping $\tvh\to\tv$ is an isomorphism of posets. }\e  

\begin{proof}
By Theorem \ref{tvBij}, the mapping is bijective, and it obviously preserves the ordering. Thus, it suffices to prove that the inverse mappping preserves the ordering too. This is an immediate consequence of Theorem A and Lemma \ref{inherite}.
\end{proof}

\end{document}